\documentclass[a4paper,11pt]{amsart}
\usepackage{amssymb,amsfonts,amsxtra,amscd}
\usepackage[all]{xy}
\usepackage{fullpage}
\theoremstyle{plain}
\newtheorem{theorem}{Theorem}[section]
\newtheorem{lemma}[theorem]{Lemma}
\newtheorem{cor}[theorem]{Corollary}
\newtheorem{prop}[theorem]{Proposition}
\theoremstyle{definition}
\newtheorem{defi}[theorem]{Definition}
\newtheorem{example}[theorem]{Example}
\theoremstyle{remark}
\newtheorem{rem}[theorem]{Remark}
\numberwithin{equation}{section}
\newcommand{\p}{\ensuremath{\mathcal H}}
\newcommand{\Z}{\ensuremath{\mathbb Z_p}}
\newcommand{\A}{\ensuremath{\mathcal H}}
\newcommand{\D}{\ensuremath{\mathcal D}}
\newcommand{\All}{\ensuremath{{\mathcal H}_{ll}}}
\newcommand{\Alr}{\ensuremath{{\mathcal H}_{lr}}}
\newcommand{\Arl}{\ensuremath{{\mathcal H}_{rl}}}
\newcommand{\Arr}{\ensuremath{{\mathcal H}_{rr}}}
\newcommand{\F}{\ensuremath{\mathbb{F}_p}}
\newcommand{\lra}{\ensuremath{\longrightarrow}}
\newcommand{\noproof}
{\begin{flushright} \ensuremath{\square} 
\end{flushright}}
\newcommand{\inlim}[2]{\ensuremath{\varprojlim_{#1} #2}}
\newcommand{\dilim}[2]{\ensuremath{\varinjlim_{#1} #2}}
\DeclareMathOperator{\Hom}{Hom}
\DeclareMathOperator{\Map}{Map}
\DeclareMathOperator{\coker}{coker}
\DeclareMathOperator{\sgn}{sgn}
\DeclareMathOperator{\hocolim}{hocolim}
\begin{document}
\begin{abstract}
We study the structure of the formal groups associated to the Morava $K$-theories of integral Eilenberg-Mac Lane spaces. The main result is that every formal group in the collection $\{K(n)^*K({\mathbb Z}, q), q=2,3,\ldots\}$  for a fixed $n$ enters in it together with its Serre dual, an analogue of a principal polarization on an abelian variety. 
We also identify the isogeny class of each of these formal groups over an algebraically closed field. These results are obtained with the help of the Dieudonn\'e correspondence between bicommutative Hopf algebras and Dieudonn\'e modules. We extend P. Goerss's results on the bilinear products of such Hopf algebras and corresponding Dieudonn\'e modules.
\end{abstract}
\title{Dieudonn\'e modules and $p$-divisible groups associated with Morava $K$-theory of Eilenberg-Mac Lane spaces}
\author{Victor Buchstaber \and Andrey Lazarev}
\address{Steklov Mathematical
Institute, Russian Academy of Sciences, Gubkina 8 Moscow 119991,
Russia.}
\email{buchstab@mendeleevo.ru}
\address{Mathematics Department,University of Bristol, Bristol, BS8 1TW, England.}
\email{a.lazarev@bristol.ac.uk}

\subjclass[2000]{}
\thanks{The authors were partially supported by the EPSRC grant GR/R84276/01.}
\maketitle

\section{Introduction}
The theory of formal groups gave rise to a powerful method for solving various problems of algebraic topology thanks to fundamental works by Novikov 
\cite{nov} and Quillen \cite{qui}. Formal groups in topology arise when one applies a complex oriented cohomology theory to the infinite complex projective space ${\mathbb C}P^\infty$. However the formal groups obtained in this way are all one-dimensional and so far the rich and intricate theory of higher dimensional formal groups remained outside of the realm of algebraic topology. One could hope to get nontrivial examples in higher dimensions by applying a generalized cohomology to an $H$-space. For most known cohomology theories and $H$-spaces this hope does not come true, however there is one notable exception. Quite surprisingly, the Morava $K$-theories  applied to integral Eilenberg-Mac Lane spaces give rise to formal groups in higher dimensions.  Moreover, these formal groups are exceptionally good in the sense that they have \emph{finite height}. 

This striking result belongs to Ravenel and Wilson \cite{ravwil} who used it to prove the so-called Conner-Floyd conjecture. However until now there has not been a systematic study of the remarkable collection of formal groups discovered by Ravenel and Wilson. This study is our main objective in this paper.

The main tool for Ravenel and Wilson was the notion of a Hopf ring and its behaviour in spectral sequences.
The definition of a Hopf ring was recently put in a conceptual framework by Goerss by introducing a suitable symmetric monoidal category for bicommutative Hopf algebras in \cite{Goe}. We make substantial use of the results of this paper. 

It is well-known that the most effective way to study formal groups, particularly those of finite height or, more generally, \emph{p-divisible groups} is via the \emph{Dieudonn\'e functor} which associates to a formal group a module over a certain  ring called the Dieudonn\'e ring cf. \cite{man}, \cite{Dem}. Goerss supplied the category of Dieudonn\'e modules with a monoidal structure and showed that the Dieudonn\'e functor is monoidal. We use this technique to study the structure of Ravenel-Wilson formal groups.

Our main result is that the spectrum multiplication \[K({\mathbb Z}/p^\nu, q)\wedge K({\mathbb Z}/p^\nu, n-q)\rightarrow K({\mathbb Z}/p^\nu, n)\] induces a kind of Poincar\'e duality on $K(n)_*K({\mathbb Z}/p^\nu, -)$ where $K(n)$ is the $n$th Morava $K$-theory. More precisely, we show that the Hopf algebras $K(n)_*K({\mathbb Z}/p^\nu, q)$ and $K(n)_*K({\mathbb Z}/p^\nu, n-q)$ are dual to each other for $n$ odd and `twisted dual' for $n$ even (precise formulations are found in the main text). Moreover, the formal groups $K(n)^*K({\mathbb Z}, q+1)$ and $K(n)^*K({\mathbb Z}, n-q+1)$ are Serre dual to each other and we identify explicitly the isogeny classes of these formal groups over an algebraic closure of $\F$, the field of $p$ elements.

The main ingredient in the proof is the theorem of Ravenel and Wilson which shows that the collection of Hopf algebras $K(n)_*K({\mathbb Z}/p^\nu, -)$ forms an \emph{exterior} Hopf ring on $K(n)_*K({\mathbb Z}/p^\nu, 1)$. The `Poincar\'e duality' mentioned above is not a formal consequence of the Ravenel-Wilson theorem, though, but follows from rather exceptional properties of the Hopf algebra $K(n)_*K({\mathbb Z}/p^\nu, 1)$.

The paper is organized as follows. In sections 2 and 3 we introduce the so-called $\boxtimes$-product (or bilinear product) in the category of bicommutative Hopf algebras $\mathcal H$. We mostly follow Goerss's paper \cite{Goe}; however our construction of the bilinear product is more explicit than his and we provide some instructive examples. In section 4 we discuss an appropriate version of the Cartier duality which is more general than the usual one in that we do not restrict our study to finite-dimensional Hopf algebras. We also answer Paul Goerss's question to explicitly describe the internal Hom functor in $\mathcal H$. Sections 5 through 8 discuss the monoidal structure on the category of Dieudonn\'e modules $\mathcal D$ as well as the Dieudonn\'e correspondence. Again, our main source is \cite{Goe}, but we also consider the internal Hom functor and duality in $\mathcal D$. In section 9 we discuss exterior Hopf rings, exterior Dieudonn\'e algebras and their relations to generalized homology of $\Omega$-spectra. Sections 10 and 11 are devoted to the study of the structure of the exterior Hopf algebra  on $K(n)_*K({\mathbb Z}/p^\nu, 1)$ and the associated Dieudonn\'e exterior algebra. The main results are formulated and proved in Section 12. 
\subsection{Notation and conventions}
We consider bialgebras over a fixed field ${\mathbf k}$ i.e. collections of data $(H,\Delta_H, m_H, \epsilon_H, i_H)$. Here $H$ is a vector space over ${\mathbf k}$, $m_H:H\otimes_{\mathbf k} H\rightarrow H$ and $\Delta_H:H\rightarrow H\otimes_{\mathbf k} H$  are associative and coassociative multiplication and comultiplication, $\epsilon_H:H\rightarrow {\mathbf k}$ and $i_H:{\mathbf k}\rightarrow H$ are the unit and the counit respectively. We will usually omit the subscript $H$ when it is clear from the context. A bialgebra having an antipode will be called a Hopf algebra.

Additionally all bialgebras and Hopf algebras will be assumed to be commutative and cocommutative. The antipode
$H\rightarrow H$ will be denoted by $[-1]_H$. Thus, for two Hopf algebras $A, B$ the set of all Hopf algebra homomorphisms $A\rightarrow B$ is an abelian group. The addition of two homomorphisms $f,g:A\rightarrow B$ is defined as the composition
\[\xymatrix{A\ar[r]^{\Delta_{A}}&A\otimes A\ar[r]^{f\otimes g}&B\otimes B\ar[r]^{m_B}&B}.\]
The zero homomorphism is the composite map
\[\xymatrix{A\ar[r]^{\epsilon_A}&{\mathbf k}\ar[r]^{i_B}&B}.\]
The additive inverse to a homomorphism $f:A\rightarrow B$ is given by precomposing $f$ with $[-1]_B:B\rightarrow B$ or, equivalently, postcomposing it with $[-1]_A:A\rightarrow A$.  It is well-known that Hopf algebras form an abelian category which will be denoted by $\mathcal H$. The zero object in $\mathcal H$ is the ground field $\mathbf k$.

If $A=B$ then the set of endomorphisms of $A$ is also a ring with respect to the composition of endomorphisms. The ring of integers $\mathbb{Z}$ maps canonically into this ring and we denote by $[n]_A$ or simply by $[n]$ the image of $n\in \mathbb{Z}$ under this map.

For an element $x$ in a Hopf algebra we will write $\Delta(x)=\sum x^{(1)}\otimes x^{(2)}$. Similarly the $n$-fold diagonal $\Delta^n:H\rightarrow H^{\otimes n}$ will be written as $\sum x^{(1)}\otimes\ldots \otimes x^{(n)}.$ Our ground field ${\mathbf k}$ will have characteristic $p$ unless indicated otherwise; from Section \ref{twisted} onwards $p$ will be odd. In Sections 2 through 4 the symbol $\otimes$ will stand for $\otimes_{\mathbf k}$; later on all unmarked tensor products are assumed to be taken over the $p$-adic integers $\Z$. We will denote by $\F$ the field consisting of $p$ elements, and by $\mathbb Q_p$ the field of $p$-adic rational numbers.
For a Hopf algebra $H$ we will denote by $F:H\rightarrow H$ the Frobenius morphism
$x\mapsto x^p$ and by $V:H\rightarrow H$ the Verschiebung.
\section{Bilinear products of Hopf algebras}
In this section we will introduce the operation in the category of Hopf algebras which models the tensor product of abelian groups.  
 \begin{defi}
Let $H_1,H_2,K$ be Hopf algebras. Let $\phi$ be a morphism of coalgebras
\[\phi:H_1\otimes H_2\rightarrow K.\]
We will write $x\circ y$ for $\phi(x,y)$. Then $\phi$ is called a \emph{bilinear map}
if the following axioms hold for all $x,y\in H_1, z,w\in H_2$: \begin{enumerate}\item
$xy\circ z=\sum(x\circ z^{(1)})(y\circ z^{(2)})$\item
$x\circ wz=\sum(x^{(1)}\circ w)(x^{(2)}\circ z)$\item
$x\circ 1=\epsilon(x)\cdot 1$\item
$1\circ z=\epsilon(z)\cdot 1$\end{enumerate}
 \end{defi}
The bilinear product of two Hopf algebras is defined with the help of a suitable universal property with respect to bilinear maps. More precisely:
\begin{defi}\label{box}
For two Hopf algebras $H_1, H_2$ their \emph{bilinear product} $H_1\boxtimes H_2$ is the unique Hopf algebra together with a bilinear map
\[\gamma:H_1\otimes H_2\rightarrow H_1\boxtimes H_2\]
such that for any bilinear map $H_1\otimes H_2\rightarrow K$ there exists a unique Hopf algebra map $H_1\boxtimes H_2\rightarrow K$ making commutative the following diagram
\[\xymatrix{H_1\otimes H_2\ar[r]^\gamma\ar[d]&H_1\boxtimes H_2\ar[dl]\\K}\]
\end{defi}
Of course, one still needs to prove that $H_1\boxtimes H_2$ satisfying the afore-mentioned universal property exists. A proof of this result is contained in \cite{Goe}. Below we will give another, hopefully more explicit description of the bilinear product.

Let $S({H_1}\otimes {H_2})$ be the symmetric algebra on ${H_1}\otimes {H_2}$. For $a_1\in {H}_1, a_2\in {H}_2$ we will write the elements $a_1\otimes a_2\in S({H_1}\otimes {H_2})$ as $a_1\circ a_2$. 
The algebra $S({H_1}\otimes{H_2})$ has a unique structure of a bialgebra defined by the requirement that the 
canonical inclusion $H_1\otimes H_2\rightarrow S({H_1}\otimes {H_2})$ be a map of coalgebras. Explicitly, for $a_1\in H_1,a_2\in H_2$ we have
\begin{equation}\label{diag}\Delta(a_1\circ a_2)=\sum a_{1}^{(1)}\circ a_{2}^{(1)}\otimes a_{1}^{(2)}\circ a_{2}^{(2)}.\end{equation} 

Consider the ideal $J$ in $S(H_1\otimes H_2)$ generated by the elements
\begin{enumerate}\item $a_1(x,y,z)=xy\circ z-\sum(x\circ z^{(1)})(y\circ z^{(2)})$ \item
$a_2(x,w,z)= x\circ wz-\sum(x^{(1)}\circ w)(x^{(2)}\circ z)$\item
$b_1(x)=x\circ 1-\epsilon(x)\cdot 1$ \item $b_2(z)=1\circ z-\epsilon(z)\cdot 1.$ \end{enumerate}
Let us show that $J$ is a coideal with respect to the coalgebra structure in $S(H_1\otimes H_2)$. We will restrict ourselves with checking the elements $a_1(x,y,z)$ and $b_1(x)$ only; the proof for $a_2(x,w,z)$ and $b_2(z)$ is similar. We have
\begin{align*}\Delta b_1(x)=&\sum x^{(1)}\circ 1\otimes x^{(2)}\circ 1-\epsilon(x)1\otimes 1\\
=&\sum(b_1(x^{(1)})+\epsilon(x^{(1)})1)\otimes (b_1(x^{(2)})+\epsilon(x^{(2)})1)-\epsilon(x)1\otimes 1\\=& \sum b_1(x^{(1)})\otimes b_1(x^{(2)})+b_1(x^{(1)})\otimes\epsilon(x^{(2)})1+\epsilon(x^{(1)}1)\otimes b_1(x^{(2)}).\end{align*} We see that $\Delta b_1(x)\in J\otimes S+S\otimes J$. Further,
\[\Delta(xy\circ z)=(\Delta x\Delta y)\circ \Delta z=\sum(x^{(1)}y^{(1)})\circ z^{(1)}\otimes(x^{(2)}y^{(2)})\circ z^{(2)}.\] On the other hand,
\begin{align*}\Delta\sum(x\circ z^{(1)})(y\circ z^{(2)})&=\sum\Delta(x\circ z^{(1)})\Delta(y\circ z^{(2)})\\
 &=\sum[(x^{(1)}\circ z^{(1)(1)})
\otimes(x^{(2)}\circ z^{(1)(2)})] 
[(y^{(1)}\circ z^{(2)(1)})
\otimes (y^{(2)}\circ z^{(2)(2)})]\\
&= \sum(x^{(1)}\circ z^{(1)(1)})
(y^{(1)}\circ z^{(2)(1)})\otimes 
(x^{(2)}\circ z^{(1)(2)})
(y^{(2)}\circ z^{(2)(2)})
\end{align*}
Using cocommutativity of the comultiplication on $H_2$ we can rewrite the last expression as follows:
\[  \sum(x^{(1)}\circ z^{(1)(1)})
(y^{(1)}\circ z^{(1)(2)})\otimes 
(x^{(2)}\circ z^{(2)(1)})
(y^{(2)}\circ z^{(2)(2)})\equiv \Delta(xy\circ z)\mod (S\otimes J+J\otimes S).\]
We have the following result.
\begin{prop}
The algebra $S(H_1\otimes H_2)/J$ has the structure of a Hopf algebra where the diagonal is given by the formula (\ref{diag}). The antipode is given by the formula
\[[-1](x\circ y)=([-1]x)\circ y.\] 
\end{prop}
\begin{proof}
It only remains to check the formula for the antipode. We have
\[\xymatrix{x\circ y\ar^-{\Delta}[r]
&\sum x^{(1)}\circ y^{(1)}\otimes x^{(2)}\circ y^{(2)}
\ar^-{[-1]\otimes[1]}[r]& 
\sum([-1]x^{(1)})\circ y^{(1)}\otimes x^{(2)}\circ y^{(2)}}\]
\[=\sum[([-1]x^{(1)})x^{(2)}]\circ y=\sum \epsilon(x)1\circ y=\epsilon(x)\epsilon(y)1  \]
as required.
\end{proof} 
\begin{rem}
Similarly we can show that $[-1](x\circ y)=x\circ([-1]y)$. Furthermore, using the equality $[-1]^2=[1]=id$ it is easy to see that the following identity holds for any $x\in H_1,y\in H_2$:
\[x\circ y=[-1]x\circ[-1]y.\]
\end{rem}
Finally we have the following result which is a direct consequence of the above constructions.
\begin{cor}
The Hopf algebra $S(H_1\otimes H_2)/J$ satisfies the universal property of Definition \ref{box} and thus realizes
the bilinear product  $H_1\boxtimes H_2$. 
\end{cor}
\begin{example}\label{1}
Let $G_1, G_2$ be abelian groups. Then we have the following isomorphism of Hopf algebras:
\[{\mathbf k}[G_1]\boxtimes {\mathbf k}[G_2]\cong {\mathbf k}[G_1\otimes_{\mathbb{Z}}G_2].\] 
\end{example}
\begin{example}\label{2}For any Hopf algebra $H$ we have natural isomorphisms:
\[ {\mathbf k}[\mathbb{Z}]\boxtimes H\cong H\cong H\boxtimes {\mathbf k}[\mathbb{Z}].\] The canonical bilinear map ${\mathbf k}[\mathbb{Z}]\otimes H\rightarrow H$ is constructed as follows. Let $t^n\in {\mathbf k}[\mathbb{Z}]={\mathbf k}[t,t^{-1}]$ and $a\in H$. Then $t^n\circ a\mapsto[n](a)$. This determines a homomorphism of Hopf algebras  ${\mathbf k}[\mathbb{Z}]\boxtimes H\rightarrow H$. The inverse map is specified by $a\mapsto t\circ a$. The second isomorphism is constructed similarly.   
\end{example}
\begin{example} \label{3}
Let $H_1={\mathbf k}[x_1,\ldots,x_n],H_2={\mathbf k}[y_1,\ldots,y_k]$ where $x_i$ and $y_i$ are primitive. Then $H_1\boxtimes H_2$ is isomorphic to the polynomial algebra on primitive generators $x_i\circ y_j i=1,2,\ldots n,j=1,2,\ldots k $.
\end{example}
\begin{example}\label{4} Let $H={\mathbf k}[x]/x^p$. Then $H\boxtimes H\cong {\mathbf k}[x]$ where $x$ is primitive.
\end{example}
\begin{rem}
The isomorphisms of Examples \ref{3} and \ref{4} could be obtained directly using our explicit construction of $H_1\boxtimes H_2$. It is simpler, however, to do this using the Dieudonn\'e correspondence which will be discussed later on.
\end{rem}
\begin{rem}
We are only interested in the case when the field ${\mathbf k}$ has characteristic $p\neq 0$. Note, however, that the construction of the bilinear product goes through also in the characteristic zero case. The isomorphisms of examples  \ref{1}, \ref{2} and \ref{3} continue to hold.  Example \ref{4} has no analogue in characteristic zero. 
\end{rem}
Summing up the above discussion we have the following result.
\begin{theorem}\label{monoidal1}
The category $\mathcal H$ together with the bilinear product $\boxtimes$ is a symmetric monoidal category. The unit for this monoidal structure is the Hopf algebra ${\mathbf k}[\mathbb{Z}]$.
\end{theorem}
\begin{proof}
The existence of the natural commutativity and associativity isomorphisms follows from the corresponding 
properties of the usual tensor product together with the universality of $\boxtimes$. The same arguments 
give rise to the commutativity of the hexagon and pentagon diagrams. We will refer to \cite{ML} for the definition and basic properties of symmetric monoidal categories.
\end{proof}

\section{Further properties of the bilinear product}   

We list a few basic properties of the product $\boxtimes$.  
\begin{prop}\label{basicidentity}
Let $x\in H_1,y\in H_2$. Then the following formulas hold in $H_1\boxtimes H_2$:
\begin{enumerate}\item 
$V(x\circ y)=Vx\circ Vy$
 \item $F(Vx\circ y)=x\circ Fy$\item
$F(x\circ Vy)=Fx\circ y.$\item
$([n]x)\circ y=[n](x\circ y)=x\circ[n]y$
\end{enumerate}
\end{prop}
\begin{proof} The stated formulas follow directly from the definition of the product $\boxtimes$.
Equations (1), (2) and (3) are essentially proved in \cite{ravwil}, Lemma 7.1. For $(4)$ we have the following identities:
\begin{align*}
([n]x)\circ y=&\sum (x^{(1)}\ldots x^{(n)})\circ y\\=&
\sum (x^{(1)}\circ y^{(1)})\ldots (x^{(n)}\circ y^{(n)})\\=&
[n](x\circ y)\\=&
\sum x\circ (y^{(1)}\ldots y^{(n)})\\=&
x\circ ([n]y).
\end{align*}

\end{proof}
 
Furthermore, for any two Hopf algebras $H_1,H_2$ there exists an internal Hom object $\underline\Hom(H_1,H_2)$
so that for any $K$ in $\mathcal H$ \[\Hom_\A({H_1\boxtimes K,H_2})\cong \Hom_\A(H_1,\underline\Hom(K,H_2)).\] 

We will denote by $\A_\nu$ the subcategory of $\A$ consisting of the Hopf algebras $H$ such that   $[p^\nu](a)=\epsilon(a)$ for any $a\in H$ (the equivalent condition is that $[p^\nu]$ is a zero element in the ring $\Hom_\A(H,H)$. Using the language of algebraic geometry one can say that $H$ represents a $p^\nu$-torsion group scheme. 
The union of all subcategories $\A_\nu$ will be denoted by $\A_\infty$.

Similarly we denote by $\A(N)$ the subcategory of $\A$ formed by those Hopf algebras for which the $N$th iteration of the Verschiebung $V^N$ is the zero endomorhism. The union of all $\A(N)$ will be denoted by $\A(\infty)$. We will call the objects in $\A(\infty)$ \emph{irreducible Hopf algebras} (in algebraic geometry the corresponding objects are called unipotent group schemes.)

Note that the category $\A_\nu$ as well as $\A_\infty$ forms an ideal inside the symmetric monoidal category $\A$ in the sense that if $H\in \A_\nu$ then for any $A\in \A$ the Hopf algebras $H\boxtimes A$ and $A\boxtimes H$ belong to $\A_\nu$. This follows from formula (4) of Proposition  \ref{basicidentity}. 

 We also have the corresponding  statement for the internal Hom functor:
\begin{prop}
Let $H$ be a Hopf algebra in $\A_\nu$. Then for any $A\in\A$ the Hopf algebra $\underline{\Hom}(A,H)$ as well as  $\underline{\Hom}(H,A)$ also belong to $\A_\nu$.
\end{prop}
\begin{proof}
Let us prove that $\underline{\Hom}(A,H)\in \A_\nu$, the remaining case is treated similarly. Consider $\underline{\Hom}(A,H)$ as a functor of the second argument; we claim that for any $n$ the map $[n]:\underline{\Hom}(A,H)\lra \underline{\Hom}(A,H)$ is induced by the map $[n]:H\lra H$. The claim obviously implies the statement of the proposition.

We have two maps: $[n]$ and $\underline{\Hom}(A,[n]):\underline{\Hom}(A,H)\lra \underline{\Hom}(A,H)$;
to show that they coincide it suffices to check that the induced two maps
\[\phi,\psi:\Hom(C,\underline{\Hom}(A,H))\lra \Hom(C,\underline{\Hom}(A,H))\]
coincide for any $C\in \A$. Identifying  $\Hom(C,\underline{\Hom}(A,H))$ with $\Hom(C\boxtimes A,H)$ we see that
$\phi$ and $\psi$ are both equal to the map $C\lra \underline{\Hom}(A,H)$ which is induced by $[n]:C\lra C$.
\end{proof}
The category $\A_\nu$ (but not $\A(N)$) has a unit making it a symmetric monoidal category:
\begin{prop}
The unit in the monoidal category $\A_\nu$ is the Hopf algebra ${\mathbf k}[\mathbb{Z}/p^\nu]$.
\end{prop}
\begin{proof}
Let $t\in {\mathbf k}[\mathbb{Z}/p^\nu]$ be the generator in $\mathbb{Z}/p^\nu$.  The unit map ${\mathbf k}[\mathbb{Z}/p^\nu]\boxtimes H\rightarrow H$ is specified by $t^n\circ a\mapsto[n]a$ for $a\in H$. This is a well-defined map since
\[\epsilon(a)=[p^\nu]a=t^{p^\nu}\circ a=1\circ a.\] The inverse map is defined as $a\mapsto t\circ a$.
\end{proof}
We will finish this section by introducing another  subcategory \A$^f$ inside {\A} consisting of \emph{finite dimensional} Hopf algebras. In the algebro-geometric literature they are known under the name of \emph{finite abelian group schemes}. It is well-known (Demazure-Gabriel, \cite{GD}) that the category \A$^f$ splits as a direct product:
\[\A^f=\All\times \Alr\times\Arl\times \Arr.\]
Here {\All} consists of Hopf algebras which are local with local dual, {\Alr} are the Hopf algebras which are local and whose duals are reduced (have no nilpotent elements), {\Arl} are the Hopf algebras which are reduced with local duals and {\Arr} stands for those Hopf algebras which are reduced together with their duals.
Let us give some typical examples of Hopf algebras in each of the four categories listed above. 
\begin{enumerate} \item $\A_{ll}$:
${\mathbf k}[x]/x^p$ where $x$ is primitive;\item $\A_{lr}$: a group algebra of a finite abelian $p$-group:\item $\A_{rl}$: a ${\mathbf k}$-dual to the group algebra of a finite abelian $p$-group;\item $\A_{rr}$: a group algebra of a finite abelian group whose order is coprime to $p$. 
\end{enumerate}
\begin{rem}Note that $\All\times \Alr\times\Arl$ coincides with $\A_\infty^f:=\A^f\bigcap\p$. The reason for considering the category $\A_\infty^f$ is that it is closed with respect to the Cartier duality which will be discussed in the next section and that it behaves well with respect to the Dieudonn\'e correspondence. However observe, that $\A^f$ or $\A_\infty^f$ are not closed with respect to $\boxtimes$.  Consider, e.g. the Hopf algebra $H={\mathbf k}[x]/x^p$ with $\Delta(x)=x\otimes 1+1\otimes x$. Then $H\boxtimes H\cong {\mathbf k}[x]$ according to Example \ref{4}, in particular, $H\boxtimes H$ is not finite-dimensional. 
\end{rem}

\section{Duality for Hopf algebras}
We will now describe a version of the Cartier duality in the category $\A$ and some of its subcategories.  We start by defining for any Hopf algebra $H\in \mathcal H$ its dual Hopf algebra $H^0$. This construction makes sense for an arbitrary (not necessarily commutative or cocommutative) Hopf algebra over a field. 

Note that if $H$ is not finite dimensional over ${\mathbf k}$ then $H^*=\Hom(H,{\mathbf k})$ is a \emph{topological} Hopf algebra
rather than a usual Hopf algebra. That means that the comultiplication $\Delta:H^*\lra H^*\hat{\otimes}H^*:=(H\otimes H)^*$ does not take its values in $H^*\otimes H^*$ as required for a Hopf algebra but in a bigger space $H^*\hat{\otimes}H^*$. We say that $A\subset H^*$ is a Hopf subalgebra of $H^*$ if $A$ is a subalgebra and $\Delta(A)\in A\otimes A\subset H^*\hat{\otimes}H^*$. In addition, we require that $A$ be closed under the antipode. Clearly then, $A$ itself is a Hopf algebra.

\begin{defi}
Let $H\in \mathcal H$. Define the (Cartier) dual Hopf algebra $H^0$ to be the union of all Hopf subalgebras inside $H^*$.
\end{defi} 
\begin{rem}
Note that if $A,B$ are Hopf subalgebras in $H^*$ then so is $A\cdot B$, the set of all linear combinations of products of elements in $A$ and $B$.  It implies that the union of all Hopf subalgebras in $H^*$ is again a Hopf algebra.

Clearly, if $H$ is finite-dimensional then $H^0$ coincides with $H^*$.  However, in general $H^0$ could be very complicated. Consider, e.g. the Hopf algebra ${\mathbf k}[t]$ where the generator $t$ is primitive. It is well-known that the graded dual to ${\mathbf k}[t]$ (where $t$ is taken to be a homogeneous element of degree $2$) is $\Gamma$, the algebra of divided powers. However, ${\mathbf k}[t]^0$ is much bigger than $\Gamma$, in particular it always contains ${\mathbf k}[{\mathbf k}]$, the ${\mathbf k}$-group ring of ${\mathbf k}$ considered as an additive group. To see that observe that there is a natural evaluation map ${\mathbf k}[t]\lra \Map({\mathbf k},{\mathbf k})$ where $\Map({\mathbf k},{\mathbf k})$ stands for the set of all maps of sets ${\mathbf k}\lra {\mathbf k}$. It is clear that $\Map({\mathbf k},{\mathbf k})$ is a topological Hopf algebra which is dual to ${\mathbf k}[{\mathbf k}]$. Taking the continuous dual to the above map we arrive at the inclusion ${\mathbf k}[{\mathbf k}]\lra H^*$. Related questions are discussed in authors' paper \cite{BL}.
\end{rem}
\begin{defi}
The dualizing object $D_h$ in $\A$ is defined as $D_h:={\mathbf k}[\mathbb{Z}]^0$, the dual Hopf algebra to the unit object ${\mathbf k}[\mathbb{Z}]$.
 Similarly define the dualizing object in $\p_\nu$ to be $D(\nu)_h:={\mathbf k}[\mathbb{Z}/p^\nu]^*$, the ${\mathbf k}$-linear dual to the unit object in $\p_\nu$.
\end{defi}
\begin{rem}
The subscript $h$ in the notation for the dualizing object stands for `Hopf' and its purpose is to distinguish it from the dualizing object in the category of Dieudonn\'e modules. We will suppress this subscript as well as the dependence on $\nu$ in cases when no confusion is possible.
\end{rem}  
\begin{defi} 
Let $A, B$ be Hopf algebras in $\A$. A map $\langle,\rangle:A\otimes B\rightarrow {\mathbf k}$ is called a \emph{bilinear pairing} if the following axioms hold for any $a, a_1, a_2\in A, b,b_1,b_2\in B$:
\begin{enumerate}
\item $\langle a, b_1b_2\rangle=\langle \Delta a,b_1\otimes b_2\rangle$\item
$\langle a_1a_2, b\rangle=\langle a_1\otimes a_2,\Delta b\rangle$\item
$\langle 1,b\rangle=\epsilon(b)$\item
$\langle a, 1\rangle=\epsilon(a) $
\end{enumerate} \end{defi}
Clearly, a bilinear pairing $A\otimes B\rightarrow {\mathbf k}$ is equivalent to a map of Hopf algebras $A\rightarrow B^0$.
In what follows we treat the Cartier  duality in the category $\p_\nu$. However Lemma \ref{pairing}, Definition \ref{dual} and Theorem \ref{Car} have obvious analogues, with similar proofs in the category \A.
 \begin{lemma}\label{pairing}
A bilinear pairing $A\otimes B\rightarrow {\mathbf k}$ for $A,B\in \p_\nu$ determines and is determined by a map of Hopf algebras
$A\boxtimes B\rightarrow D(\nu)$.
\end{lemma}
\begin{proof}
Let $\phi: A\boxtimes B\rightarrow D(\nu)$ be a Hopf algebra map. Taking its linear dual we obtain a map of algebras $D(\nu)^*={\mathbf k}[\mathbb{Z}/p^\nu]\rightarrow (A\boxtimes B)^*$ (since $\phi$ is a map of coalgebras). The latter gives rise to an invertible element in $(A\boxtimes B)^*$ whose $p^\nu$th power is $1$. We can consider this element as a map $A\otimes B\rightarrow {\mathbf k}$. The axioms 1-4 for the bilinear pairing follow from the corresponding axioms for a bilinear map.

Conversely, a bilinear pairing $\langle,\rangle:A\otimes B\rightarrow {\mathbf k}$ could be considered as an element $f\in(A\otimes B)^*$.  We need to show that $f^{p^\nu}=1$ in $(A\otimes B)^*$. Let $a\otimes b\in A\otimes B$. We have:
\begin{align*}
f^{p^\nu}(a\otimes b)&=\sum f(a^{(1)}\otimes b^{(1)})\ldots f(a^{(p^\nu)}\otimes b^{(p^\nu)})\\
&=\sum \langle a^{(1)},b^{(1)}\rangle \ldots \langle a^{(p^\nu)}, b^{(p^\nu)}\rangle\\
&=\langle a^{(1)}\ldots a^{(p^\nu)},  b\rangle\\
&=\langle[p^{\nu}](a), b\rangle \\
&=\langle\epsilon(a), b\rangle\\
&=\epsilon(a)\epsilon(b).
\end{align*}
as required.
\end{proof}

 \begin{defi}\label{dual}
For $H\in \p_\nu$ define its dual $D_hH$ as $D_hH:=\underline\Hom(H,D(\nu))$.
\end{defi}
Following our customary abuse of notation we will usually shorten $D_hH$ to $DH$.
Then we have the following result.
\begin{theorem}\label{Car}
There is a natural isomorphism of Hopf algebras $DH\cong H^0$ in $\p_\nu$.
\end{theorem}
 \begin{proof}
Let $A$ be a Hopf algebra in $\p$. We have a natural isomorphism
$\Hom_{\p_\nu}(A,DH)\cong \Hom_{\p_\nu}(A\boxtimes H, D)$. By the previous lemma the set  $\Hom_{\p_\nu}(A\boxtimes H,D)$
is in natural $1-1$ correspondence with the set of bilinear pairings $A\otimes H\rightarrow {\mathbf k}$ and the latter bijects with the set of topological Hopf algebra maps $A\rightarrow H^*$. Thus, there is a natural bijection of sets
\[\Hom_{\p_\nu}(A,DH)\approx\Hom_{\p_\nu}(A,H^0)\]
and therefore $DH\cong H^0$.
\end{proof}
\begin{cor}
The Cartier duality functor restricts to the subcategory $\p^f$. \begin{proof}Indeed, it interchanges the categories 
$\Alr$ and $\Arl$ and maps the category $\All$ to itself. \end{proof}
\end{cor}
\begin{cor}
For any finite-dimensional Hopf algebra $H$ in $\p_\nu$ there is a natural isomorphism $H\cong DDH$.
\end{cor}
\noproof
\begin{cor}
Let $H_1,H_2$ be Hopf algebras in $\p_\nu$ and assume that $H_2$ is finite dimensional. Then  \[\underline{\Hom}(H_1,H_2)\cong (H_1\boxtimes H_2^*)^0.\]
\end{cor}
\begin{proof}
We have the following sequence of isomorphisms of Hopf algebras
\begin{align*}
\underline{\Hom}(H_1,H_2)&\cong\underline{\Hom}(H_1,DDH_2)\\&\cong
\underline{\Hom}(H_1\boxtimes DH_2, D)\\&\cong D(H_1\boxtimes DH_2)\\&\cong (H_1\boxtimes H_2^*)^0
\end{align*}
\end{proof}
The above formula for the internal Hom functor answers the question of Paul Goerss, cf. \cite{Goe}, section 5, at least for finite dimensional $p$-torsion Hopf algebras. In the next section we will give a corresponding formula for the internal Hom of Dieudonn\'e modules.

\section{The bilinear product of Dieudonn\'e modules}
In this section we specialize ${\mathbf k}=\F$. There is little doubt that all our constructions could be generalized to 
the case of an arbitrary perfect field, however our topological examples do not require this level of generality and we restrict ourselves with considering the case of a prime field only.
\begin{defi}
The category $\D$ is the category of modules over the Dieudonn\'e ring $R=\Z[V,F]/(VF-p)$. The objects in $\D$ will be called Dieudonn\'e modules. The subcategory 
$\D(N)$ consists of those Dieudonn\'e modules for which $V^N$ acts trivially. The subcategory $\D_\nu$ of $\D$ consists
of those Dieudonn\'e modules for which $p^\nu$ acts trivially. The union of $\D(N)$ will be denoted by $\D(\infty)$.
\end{defi}
Clearly $\D, \D(N),\D_\nu$ are abelian categories. We will next introduce the notion of a \emph{bilinear map} in $\D$ similar to the bilinear map of Hopf algebras discussed in the second section of the paper.
\begin{defi}\label{bilinearD}
Let $M, N, L$ be $R$-modules. A map $f:M\otimes N\lra L$ is called a bilinear map if it is $\mathbb{Z}_p$-bilinear and
\begin{enumerate}\item
$Ff(Vm\otimes n)=f(m\otimes Fn)$\item
$Ff(m\otimes Vn)=f(Fm\otimes n)$\item
$Vf(m\otimes n)=f(Vm\otimes Vn)$.
\end{enumerate} 
\end{defi}
Just as before, the notion of a bilinear map leads naturally to the notion of a \emph{bilinear product} in the category of Dieudonn\'e modules.
\begin{defi}
For two Dieudonn\'e modules $M$ and $N$ we define their bilinear product $M\boxtimes N$ as the unique Dieudonn\'e module supplied with a bilinear map $M\otimes N\rightarrow M\boxtimes N$ such that for any 
 bilinear map $M\otimes N\rightarrow L$ there exists a unique $R$-module map $M\boxtimes N\rightarrow L$ making commutative the following diagram
\[\xymatrix{M\otimes N\ar[r]\ar[d]&M\boxtimes N\ar[dl]\\L}.\]\end{defi}
An explicit description of $M\boxtimes N$ is given in the following construction.   

Let $M,N\in\D$ and consider $M\otimes N$ as an $\Z[V]$-module so that 
\[V(x\otimes y)=Vx\otimes Vy.\]  Then $R\otimes_{\Z[V]}(M\otimes N)$ has an obvious structure of a left $R$-module. Set
\[M\boxtimes N:=R\otimes_{\Z[V]}(M\otimes N)/\sim\]
where $\sim$ is the $R$-submodule generated by the following elements:
\begin{equation}\label{rel1}
F\otimes Vm\otimes n-1\otimes m\otimes Fn;\end{equation}
\begin{equation}\label{rel2}
 F\otimes m\otimes Vn-1\otimes Fm\otimes n.\end{equation}
Hence the structure of an $R$-module on $M\boxtimes N$ as follows:
\[F(r\otimes m\otimes n)=Fr\otimes m\otimes n;\]
\[V(r\otimes m\otimes n)=r\otimes Vm\otimes Vn.\]
The element $1\otimes m\otimes n\in M\boxtimes N$ will be denoted by $m\circ n$.

The following result is an analogue of Theorem \ref{monoidal1}. Its proof is a direct inspection of definitions.
\begin{theorem}
The category $\D$ together with the bilinear product $\boxtimes$ is a symmetric monoidal category. The unit for this monoidal structure is $R$-module ${\mathbb I}_\D=\Z$ where $F$ acts as a multiplication
by $p$ and $V$ is the identity automorphism.
\end{theorem}\noproof

The product $\boxtimes$ also determines a monoidal structure in $\D_\nu$ and $\D(N)$. The unit in $\D_\nu$ is the module
${\mathbb I}_{\D(\nu)}=\Z/p^\nu$ where, as above, $F$ acts as a multiplication
by $p$ and $V$ is the identity automorphism. The category $\D(N)$ has no unit.
\begin{rem}
The structure of the bilinear product of two $R$-modules is not so obvious even in the simplest cases. For example one can prove that $R\boxtimes R$ is isomorphic to the direct sum of copies of $R$ with generators $1\circ 1, V^l\circ 1, 1\circ V^l$ where $l=1,2,\ldots$.
\end{rem}
\section{Duality for Dieudonn\'e modules}

We will start our treatment of duality in the category $\D$ by introducing the internal Hom functor. 
\begin{defi}\label{internal}Let $M,N\in \D$ or $M, N\in \D_\nu$. Define $\underline{\Hom}(M,N)$ to be the subgroup of $\Hom(R\otimes M,N)$
consisting of such $f:R\otimes M\lra N$ for which\begin{enumerate}\item
$Ff(Vr\otimes m)=f(r\otimes Fm)$\item
$Ff(r\otimes Vm)=f(Fr\otimes m)$\item
$Vf(r\otimes m)=f(Vr\otimes Vm)$\end{enumerate}
Furthermore, define the structure of an $R$-module on  $\underline{\Hom}(M,N)$ by the following formula:
\begin{eqnarray}\label{R-mod}(r\cdot f)(a\otimes m)=f(ra\otimes m).\end{eqnarray} 
\end{defi}
\begin{rem}
It is clear that $\underline{\Hom}(M,N)$ is isomorphic to $\Hom_R(R\boxtimes M,N)$ as an abelian group. It implies, in particular, that $\underline{\Hom}(M,N)$ is a $\Z$-module. However the structure of an $R$-module on $\underline{\Hom}(M,N)$ is different from that on $\Hom_R(R\boxtimes M,N)$ 
\end{rem}
\begin{prop}
Formulas (\ref{R-mod}) determine the structure of an $R$-module on $\underline{\Hom}(M,N)$.
\end{prop}
\begin{proof}
One needs to check that $(F\cdot f)(a\otimes m)$ and $(V\cdot f)(a\otimes m)$ belong to $\underline{\Hom}(M,N)$
i.e. that they satisfy the formulas (1)-(3) of Definition \ref{internal} and that the relations $VF=FV=p$ hold.
This verification is completely straightforward. 
\end{proof}
\begin{prop}
There is a natural isomorphism in \D:
\begin{equation}\label{conj}\Hom_R(M\boxtimes N,L)\cong \Hom_R(M,\underline{\Hom}(N,L)).\end{equation}
\end{prop}
\begin{proof}
The $R$-module $\underline{\Hom}(N,L)$ is, by definition, a subgroup in $\Hom(R\otimes N,L)$. Therefore $\Hom_R(M,\underline{\Hom}(N,L))$ is identified with a certain subgroup inside
\begin{align*}\Hom_R(M,\Hom(R\otimes N,L))&\cong \Hom(M\otimes_R(R\otimes N),L)\\
&\cong \Hom(M\otimes N,L).\end{align*}
More precisely, an easy inspection shows that this subgroup consists of $f\in \Hom(M\otimes N,L)$ for which
the identities (1)-(3) of Definition \ref{bilinearD} hold.
It follows that the collection of such $f$ is isomorphic to $\Hom_R(M\boxtimes N,L)$.

So we showed that (\ref{conj}) is an isomorphism of abelian groups. Let us now check that this is an $R$-module isomorphism. Let $f\in \Hom_R(M\boxtimes N,L)$. We will consider $f$ as a map $M\otimes N\lra L$. Denote by $\bar{f}$ the corresponding
homomorphism 
\[M\lra\underline{\Hom}(N,L)\subset \Hom(R\otimes N,L).\] 
Thus, $\bar{f}(m)(r\otimes n)=f(rm\otimes n).$ Let $a\in R$. We have:
\begin{align*}
\bar{f}(am)(r\otimes n)&=f(ram\otimes n)\\
&=\bar{f}(m)(ar\otimes n)\\&=(a\cdot\bar{f})(m)(r\otimes n).
\end{align*}
\end{proof}

Finally, we introduce the notion of the Cartier duality in the category $\D_\nu$.
\begin{defi}
The dualizing module in the category  $\D_\nu$ is the group 
$D_d(\nu):={\mathbb Z}/p^\nu{\mathbb Z}$ where $V$ acts as multiplication by $p$ and $F$ is the 
identity automorphism. 
For $M\in \D_\nu$ define its dual Dieudonn\'e module $D_dM$ as $\underline{\Hom}(M,D_d(\nu))$.\end{defi}
\begin{rem}
The subscript $d$ in the definition of the dualizing module is supposed to distinguish it from the dualizing object for Hopf algebras. We will suppress this subscript as well as the dependence on $\nu$ whenever practical.
\end{rem}

We now introduce the notion of a \emph{bilinear pairing} of Dieudonn\'e modules  analogous to the corresponding notion in the category of Hopf algebras.
\begin{defi}\label{pairD}
A bilinear pairing of two Dieudonn\'e modules $M,N\in \D_\nu$ is a map
\[\{,\}:M\otimes N\rightarrow \mathbb{Z}/p^\nu\]
such that for any $m\in M$ and $n\in N$ \begin{enumerate}
\item $\{m,Fn\}=\{Vm,n\}$\item
 $\{Fm,n\}=\{m,Vn\}$.
\end{enumerate}
\end{defi}
\begin{rem}
Obviously a bilinear pairing $M\otimes N\rightarrow {\mathbb Z}/p^\nu$ nothing but a bilinear map 
$M\otimes N\rightarrow D$ or, equivalently,
a map of Dieudonn\'e modules $M\boxtimes N\rightarrow D$.
\end{rem}

 We have the following result whose proof is a simple check.
\begin{prop}\label{dua}
For $M\in \D_\nu$ the module $DM$ can be identified with the abelian group 
$\Hom(M,\mathbb Z/p^\nu)$.
The action of the operators $F$ and $V$ is specified by the formulas
\[(F\cdot f)(m)=f(Vm);\]
\[(V\cdot f)m=f(Fm).\]
\end{prop}\noproof
\begin{rem}
It is also possible to consider the Cartier duality in the whole category $\D$. In this case the dualizing module $D$ will be the abelian group $\mathbb Q_p/\Z$ where $V$ and $F$ act as before. The analogue of  Proposition \ref{dua} continues to hold in this context. 
\end{rem}

\section{The Dieudonn\'e correspondence}
Recall that the Witt Hopf algebra $W_n^{\mathbb Z}=\mathbb{Z}[x_0,\ldots,x_n]$ admits a unique Hopf algebra
structure over $\mathbb{Z}$ for which the Witt polynomials $P_k=x_0^{p^k}+px_1^{p^{k-1}}+\ldots+p^kx_k$ are primitive. We will write $W_n$ for $W_n^{\mathbb Z}\otimes \F$.

There is a map of Hopf algebras $\bar{V}_n:W_{n+1}\lra W_n$ defined as $\bar{V}_n(x_0)=0$ and $\bar{V}_n(x_i)=x_{i-1}$ for $i>0$. Then the classical Dieudonn\'e theorem cf. \cite{GD} states:
\begin{theorem}
 The functor ${\mathcal F}:\A\lra\D$:
 \[H\mapsto \dilim{n}\Hom_{\A}(W_n,H)\]
 establishes an equivalence of the subcategory $\A(\infty)$ of $\A$ and the subcategory $\D(\infty)$ of $\D$.
\end{theorem}
The proof of the above theorem uses the fact that $W_n$ is a projective generator of the abelian category
$\mathcal{H}(n)$.

Next recall that the categories $\A(\infty)$ and $\D$ are monoidal. The following result shows that ${\mathcal F}$ is a monoidal functor.
\begin{theorem}
There is a natural isomorphism
\[{\mathcal F}(H_1\boxtimes H_2)\cong {\mathcal F}(H_1)\boxtimes {\mathcal F}(H_2).\]
\end{theorem}
\begin{proof}
A detailed proof in the graded case is contained in Goerss's paper \cite{Goe}. Goerss's scheme carries over to the ungraded case and we will show briefly how this is done.

The first step is to construct a homomorphism of $R$-modules
\begin{equation}\label{natural}\phi:{\mathcal F}(H_1)\boxtimes {\mathcal F}(H_2)\lra {\mathcal F}(H_1\boxtimes H_2).\end{equation} Arguing as in Goerss's paper one can show that 
for a Hopf algebra $H$ over $\Z$ having a lifting of the Frobenius the Dieudonn\'e module ${\mathcal F}(H\otimes \F)$ is isomorphic to
$R\otimes _{\Z[V]}Q(H)$ where $Q(H)$ is the space of indecomposables of $H$. Next one proves that for such Hopf algebras $H_1, H_2$ there is an isomorphism 
\[Q(H_1\boxtimes H_2)\cong Q(H_1)\otimes Q(H_2).\]

Since the Hopf algebra $W_n\boxtimes W_n$ does have a lifting of the Frobenius we conclude that
 \begin{align*} {\mathcal F}(W_n\boxtimes W_n) &\cong   R\otimes _{\Z[V]}Q(W_n^{\mathbb Z}\boxtimes W_n^{\mathbb Z})\\&\cong R\otimes _{\Z[V]}Q(W_n^{\mathbb Z})\otimes Q(W_n^{\mathbb Z}).\end{align*}
Consider the element \[\kappa_n\circ\kappa_n=1\otimes x_n\otimes x_n\in R\otimes _{\Z[V]}Q(W_n^{\mathbb Z})\otimes Q(W_n^{\mathbb Z}).\]
This element represents a map of Hopf algebras 
\[\underline{\Delta}_n: W_n\lra W_n\boxtimes W_n.\]
We  will use $\underline{\Delta}_n$ to construct the  homomorphism (\ref{natural}) as follows. Without loss of generality assume that in $H_1$ and $H_2$ the operator $V^n$ is trivial. Then ${\mathcal F}(H_1)=\Hom_{\A}(W_n,H_1)$ and
${\mathcal F}(H_2)=\Hom_{\A}(W_n,H_2)$. Consider the composite map
\[\tilde{\phi}:\xymatrix{{\mathcal F}(H_1)\otimes {\mathcal F}(H_2)\cong \Hom_{\A}(W_n,H_1)\otimes \Hom_{\A}(W_n,H_2)\ar[d]\\
\Hom_{\A}(W_n\boxtimes W_n, H_1\boxtimes H_2)\ar[r]&\Hom_{\A}(W_n,H_1\boxtimes H_2)\cong {\mathcal F}(H_1\boxtimes H_2).}\]
Here the last map is induced by $\underline{\Delta}_n$. Then one can show that $\tilde{\phi}$ is a bilinear
map and therefore induces a map (\ref{natural}).

Finally, one shows that $\phi$ is an isomorphism for $H_1=W_k,H_2=W_l$ for any $k$ and $l$ and then derives that $\phi$ is an isomorphism in general.
\end{proof}
Observe that the Dieudonn\'e functor is defined on the category $\A(\infty)$  of irreducible Hopf algebras. This category is not closed with respect to the Cartier duality. There is another version of the Dieudonn\'e equivalence which is defined on the category $\A_\infty^f$ of finite dimensional $p$-torsion Hopf algebras. We will now describe this version.  Note that this category \emph{is} closed with respect to the Cartier duality.

The categories $\All$ and $\Alr$ consist of irreducible Hopf algebras and therefore the Dieudonn\'e functor is defined on them as above. We will define the Dieudonn\'e functor on $\Arl$ using the fact that $D(\Arl)=\Alr$.
Namely, for $H\in \Arl$ set ${\mathcal F}(H)=D{\mathcal F}(DH)$.

Thus, the functor ${\mathcal F}$ maps the category $\p^f_\infty$ onto the subcategory $\D^f$ of $R$-modules consisting of Dieudonn\'e module of finite length. This subcategory is the product of three subcategories $\D_{ll}, \D_{lr}$ and $\D_{rl}$. Here $\D_{ll}$ consists of those $R$-modules for which $V$ and $F$ act nilpotently, $\D_{lr}$ is the $R$-modules with nilpotent $F$ and invertible $V$ and $\D_{rl}$ is the $R$-modules with nilpotent $V$ and invertible $F$. We will sum this up in the following theorem, cf. \cite{Dem}:
\begin{theorem}\label{finite}
The functor ${\mathcal F}$ establishes an equivalence of categories $\p_\infty^f$ and $\D^f$. The functor ${\mathcal F}$ respects the Cartier
duality in $\p^f$ and $\D^f$.
\end{theorem}
\section{Twisted duality}\label{twisted}
Over the prime field $\F$ there are other candidates for a dualizing object in the categories $\mathcal H$ and $\mathcal D$ all of which become isomorphic upon passing to the algebraic closure of $\F$. In this section we consider one particular choice of a dualizing object since it will arise naturally in our study of the Morava $K$-theories of Eilenberg-Mac Lane spaces. 

From now on our ground field $\mathbf k$ will have characteristic $p>2$. We will work here in the categories ${\p}_\nu$ and ${\mathcal D}_\nu$ of Hopf algebras and Dieudonn\'e modules which are annihilated by the $\nu$th power of $p$ although one could make similar constructions in other categories.
\begin{defi}\
\begin{enumerate}\item
The twisted dualizing module $D_d^\prime$ in ${\mathcal D}_\nu$ is the abelian group ${\mathbb Z}/p^\nu$ where $F$ and $V$ act as the multiplication by $-1$ and $-p$ respectively. For a Dieudonn\'e module $M$ its twisted dual $D_d^\prime M$ is defined as $D^\prime_d M=\underline\Hom(M,D^\prime_d)$.\item
The twisted dualizing object in ${\p}_\nu$ is the Hopf algebra $D^\prime_h$ whose Dieudonn\'e module is $D_d^\prime$. For a Hopf algebra $H$ we define its twisted dual $D^\prime_hH$ as $D^\prime_hH:=\underline\Hom(M,D^\prime_h)$.
\end{enumerate} 
\end{defi}
\begin{rem}
Recall that the usual dualizing object in ${\p}_\nu$ is the dual to the group algebra of the group ${\mathbb Z}/p^\nu$. It is unlikely that one can give such a simple explicit description of $D^\prime_h$ and so we have to resort to the Dieudonn\'e correspondence instead.
\end{rem}
The notion of a twisted dualizing object leads one to introduce \emph{twisted bilinear pairing}.
\begin{defi}\label{pairTw}
A twisted bilinear pairing of two Dieudonn\'e modules $M,N\in \D_\nu$ is a map
\[[,]:M\otimes N\rightarrow \mathbb{Z}/p^\nu\]
such that for any $m\in M$ and $n\in N$ \begin{enumerate}
\item $[m,Fn]=-[Vm,n]$\item
 $[Fm,n]=-[m,Vn]$.
\end{enumerate}
\end{defi}
\begin{rem}
Obviously a twisted bilinear pairing $M\otimes N\rightarrow {\mathbb Z}/p^\nu$ is nothing but 
a bilinear map $M\otimes N\rightarrow D_d^\prime$ or, equivalently, a map of Dieudonn\'e modules $M\boxtimes N\rightarrow D_d^\prime$.
\end{rem}  

To relate two types of duality in $\p_\nu$ we need a more general form of the Dieudonn\'e correspondence, cf. \cite{GD} which we will now recall. 

Let $\mathbf k$ be an algebraic extension of $\F$ obtained by adjoining a root of some irreducible polynomial $h(x)$ so that ${\mathbf k}\cong \F[x]/(h(x))$. Denote by $\bar{h}(x)$ an integral lifting of $h(x)$. Then the ring of Witt vectors $W({\mathbf k})$ is the ring $\Z[x]/\bar{h}(x)$. It possesses a lifting of the Frobenius automorphism on $\mathbf k$ which will be denoted by $\sigma$. Then the Dieudonn\'e ring $R({\mathbf k})$ consists of all finite sums of the form
\[\sum_{i>0}\alpha_{-i}V^i+\alpha_0+\sum_{i>0}\alpha_iF^i\] where $\alpha_{i}\in W({\mathbf k})$. The multiplication law is determined by the commutation relations \[VF=FV=p;\] \[F\alpha=\sigma(\alpha)F;\] \[V\sigma(\alpha)=\alpha V\] where $\alpha\in W(\mathbf k)$. Furthermore for any Dieudonn\'e module $M$ the $R(\mathbf k)$-module $M(\mathbf k)$ is defined as $W({\mathbf k})\otimes M$ where $F$ and $V$ act according to the commutation rules above. Note that the definition of the ring $W(\mathbf k)$ can easily be extended to the case of an infinite algebraic extension by simply taking the union of $W({L})$ over subfields $L$ of $\mathbf k$ having finite degree over $\F$.

The Dieudonn\'e correspondence in this case reads as follows:
\begin{theorem}\label{general}
The category ${\mathcal H}(\infty)$ of irreducible Hopf algebras over $\mathbf k$ is equivalent to the subcategory of 
modules over $R({\mathbf k})$ for which $V$ acts nilpotently. For a Hopf algebra $H$ over $\F$ the Dieudonn\'e module corresponding to $\mathbf k\otimes H$ is isomorphic to ${\mathcal F}(H)(\mathbf k)$.
\end{theorem} 
We now have the following result.
\begin{prop}
Let $\mathbf k$ be a quadratic extension of $\F$. Then ${\mathbf k}\otimes D^\prime_h$ and ${\mathbf k}\otimes D_h$ are isomorphic Hopf algebras.
\end{prop}
\begin{proof}
The Dieudonn\'e module corresponding to ${\mathbf k}\otimes D_h$ is $W({\mathbf k})\otimes\mathbb{Z}/p^\nu$ where $F$ and $V$ act on the factor $\mathbb{Z}/p^\nu$ as the identity and multiplication by $p$ respectively.  Similarly the Dieudonn\'e module of ${\mathbf k}\otimes D_h^\prime$  is  the same abelian group $W({\mathbf k})\otimes\mathbb{Z}/p^\nu$ where now $F$ and $V$ act on the factor $\mathbb{Z}/p^\nu$ as the minus identity and multiplication by $-p$ respectively. 

Observe that since $p>2$ the extension $\mathbf k$  is obtained by adding to $\F$ the square root of a certain element $q\in \mathbb{F}_p$. Note that the Frobenius automorphism $\sigma:{\mathbf k\rightarrow \mathbf k}$ permutes the roots of the quadratic polynomial $x^2-q$, i.e. $\sigma(\sqrt{q})=-\sqrt{q}$. It follows that the map $1\otimes x\mapsto \sqrt{q}\otimes x:{\mathbf k}\otimes D_h^\prime\rightarrow {\mathbf k}\otimes D_h$ establishes an isomorphism between these two Dieudonn\'e module structures. By Theorem \ref{general} the corresponding Hopf algebras are isomorphic.
 \end{proof}
\begin{cor}
Let $H$ be a Hopf algebra in ${\p}_\nu$. Then for a quadratic extension $\mathbf k$ of $\F$ the Hopf algebras ${\mathbf k}\otimes D_hH$ and ${\mathbf k}\otimes D^\prime_hH$ are isomorphic.
\end{cor}\noproof 
  
 \section{Hopf rings, Dieudonn\'e algebras and generalized homology of $\Omega$-spectra}
In this section we introduce the notion of a Hopf ring, the corresponding notion of a Dieudonn\'e algebra and relate them to generalized homology of multiplicative 
spectra. Note that Hopf rings originally appear in the work of Milgram \cite{mil} and were later studied by Ravenel and Wilson \cite{ravwil1} in the context of
their calculation of the homology of $\Omega$-spectrum $MU$. 
\begin{defi}
A Hopf ring $A$ is a Hopf ${\mathbf k}$-algebra together with a map $\phi: A\boxtimes  A\lra A$ which is required to be 
associative: $\phi(\phi\boxtimes 1)=\phi(1\boxtimes \phi)$. A commutative Hopf ring is a Hopf ring $A$ for which
the product $\phi$ is commutative. In other words a (commutative) Hopf ring is a (commutative) monoid in the symmetric monoidal category $\A$. 
\end{defi}
\begin{rem}
Sometimes it is convenient to  require that a Hopf ring $A$ have a unit. In other words there is a map ${\mathbf k}[\mathbb{Z}]\lra A$ subject to the obvious conditions. (Of course in the category $\p_\nu$ the appropriate notion of a unit is a map ${\mathbf k}[\mathbb{Z}/p^\nu]\lra A$). 
\end{rem}
The main example of a Hopf ring is the \emph{exterior} Hopf ring on a Hopf algebra $H$ (Ravenel and Wilson use the term 'free Hopf ring').
\begin{defi}
Let $H$ be a Hopf algebra. Then its exterior Hopf ring $\Lambda_\boxtimes(H)$ is defined as
\[\Lambda_\boxtimes(H)={\mathbb I}\oplus H\oplus (H\boxtimes H)/\Sigma_2\oplus\ldots \oplus H^{\boxtimes n}/\Sigma_n\oplus\ldots.\]
Here ${\mathbb I}$ stands for the appropriate unit (${\mathbf k}[\mathbb{Z}]$ in $\A$ and  ${\mathbf k}[\mathbb{Z}/p^\nu]$ in $\p_\nu$) and $\Sigma_n$ is the symmetric group on $n$ symbols which operates on $H^{\boxtimes n}$ according  to the rule
\[\sigma(h_1\circ\ldots \circ h_n)=[-1]^{\sgn\sigma}h_{\sigma(1)}\circ\ldots \circ h_{\sigma(n)}\] where $h_1,\ldots, h_n\in H$. 

The nonunital version $\Lambda_{\boxtimes+}(H)$ is defined as
\[\Lambda_{\boxtimes+}(H)=H\oplus (H\boxtimes H)/\Sigma_2\oplus\ldots \oplus H^{\boxtimes n}/\Sigma_n\oplus\ldots.\]                                            
\end{defi}
Clearly $\Lambda_\boxtimes(H)$ and $\Lambda_{\boxtimes+}(H)$ are both Hopf rings with respect to the operation $\circ$, and $\Lambda_\boxtimes(H)$ also has a unit. The Hopf rings $\Lambda_{\boxtimes+}(H)$ and $\Lambda_\boxtimes(H)$ are not commutative, but rather skew-commutative Hopf rings in the sense that the following relation holds for any $h_1,h_2\in H$:
\[h_1\circ h_2=[-1]h_2\circ h_1.\]
The corresponding notion in the category of Dieudonn\'e modules is called the Dieudonn\'e algebra:
\begin{defi}
A (commutative) Dieudonn\'e algebra is a (commutative) monoid in the category $\D$ or $\D_\nu$.  
\end{defi}
The definition of an exterior Dieudonn\'e algebra is likewise clear:
\begin{defi}
Let $M$ be an $R$-module. Then its exterior algebra in $\D$ or $\D_\nu$ is defined as
\[\Lambda_\boxtimes(M)={\mathbb I}\oplus M\oplus (M\boxtimes M)/\Sigma_2\oplus\ldots \oplus M^{\boxtimes n}/\Sigma_n\oplus\ldots.\]
Here $\mathbb I$ is the unit in $\D$ or $\D_\nu$, i.e. the Dieudonn\'e module $\Z$ or ${\mathbb Z}/p^\nu$ with appropriate actions of $F$ and $V$.  The symmetric group $\Sigma_n$  operates on $M^{\boxtimes n}$ according  to the rule
\[\sigma(m_1\circ\ldots \circ m_n)=(-1)^{\sgn\sigma} m_{\sigma(1)}\circ\ldots \circ m_{\sigma(n)}\] where $m_1\ldots m_n\in m$.
The nonunital version $\Lambda_+(M)$ is defined as
\[\Lambda_{\boxtimes+}(M)=M\oplus (M\boxtimes M)/\Sigma_2\oplus\ldots \oplus M^{\boxtimes n}/\Sigma_n\oplus\ldots.\]
\end{defi}
Again, the exterior Dieudonn\'e algebra $\Lambda_{\boxtimes}(M)$ is not commutative, but rather skew-commutative in the obvious sense.
\begin{rem}
Let $H$ be an irreducible Hopf algebra. Then $H^{\boxtimes n}$ as well as $(H^{\boxtimes n})/\Sigma_n$ will be irreducible as well for any $n>0$. By the Dieudonn\'e correspondence the $R$-module ${\mathcal F}(\Lambda_{\boxtimes+}(H))$ is isomorphic to $\Lambda_{\boxtimes+}({\mathcal F}(H))$. In particular, it is a Dieudonn\'e algebra. If $H$ is a finite-dimensional Hopf algebra we can use the version of the Dieudonn\'e equivalence given in Theorem \ref{finite} and extend ${\mathcal F}$ to the unital 
Hopf ring $\Lambda(H)$. In this case we obtain ${\mathcal F}(\Lambda_{\boxtimes}(H))=\Lambda_{\boxtimes}({\mathcal F}(H))$.
\end{rem}
Now let $E$ be an $\Omega$-spectrum, i.e. a sequence $\{E_q$, $q=0,1,\ldots\}$ of based spaces together with 
weak equivalences $\Omega(E_{q+1})\simeq E_q$. We assume that $E$ is a ring spectrum, i.e. there exist maps of based space \begin{equation}\label{mult}E_q\wedge E_l\lra E_{q+l}\end{equation} satisfying the usual associativity axioms.

Furthermore, let $h_*(-)$ be a generalized (multiplicative) homology theory which satisfies the perfect K\"unneth formula i.e. $h_*(X\times Y)\cong h_*(X)\otimes_{h_*}h_*(Y)$ for any spaces $X$ and $Y$. Then according to \cite{ravwil1} the graded $h_*$-module $\bigoplus_nh_*(E_n)$ has the structure of a Hopf ring where the $\circ$-product
\[\circ: h_*(E_n)\otimes h_*(E_m)\lra h_*(E_{n+m})\] is induced by the mutiplication (\ref{mult}) in $E$.
Now take $E=H\mathbb{Z}/p^\nu$, the Eilenberg-Mac Lane spectrum mod $p^\nu$.  In this case $E_q=K(\mathbb{Z}/p^\nu, q)$. Next, set $h_*(-)=K(n)_*(-)$, the $nth$ Morava $K$-theory at an odd prime $p$. It will be convenient for us to consider its ungraded version $\bar{K}(n)$. It is defined as $\bar{K}(n)_*(X)=K(n)_*(X)\otimes_{K(n)_*}\F$ where $K(n)_*$ acts on $\F$ through the map $K(n)_*=\F[v_n^{\pm 1}]\rightarrow \F:v_n\mapsto 1$. Since we do not consider the graded Morava $K$-theory we will 
use the notation  $K(n)_*(-)$ for $\bar{K}(n)_*(-)$.

Now recall the following fundamental result of Wilson and Ravenel \cite{ravwil}.
\begin{theorem}
The Hopf ring ${K}(n)_*K(\mathbb{Z}/p^\nu, -)$ is the exterior Hopf ring on ${K}(n)_*K({\mathbb Z}/p^\nu,1)$.
\end{theorem}
The Hopf algebra  $H_\nu={K}(n)_*K({\mathbb Z}/p^\nu,1)$ is a well-understood object.  It is a finite dimensional Hopf algebra whose dual has the form 
$H^*_\nu=\F[[t]]/[p^\nu](t)$ where the diagonal is induced by the formal group law of the Morava $K$-theory and $[p^\nu](t)$ is the corresponding $p$-series. If one uses Hazewinkel's generators to construct the spectrum $K(n)$ then its $p$-series has the form $[p](t)=t^{p^n}$.

Note that $H_\nu$ is a finite-dimensional irreducible Hopf algebra. Let us now find its Dieudonn\'e module. Consider instead the Hopf algebra 
\[H=\inlim{\nu}H_\nu^*=K(n)^*K({\mathbb Z},2).\]
Then $H=\F[[t]]$ and $\Delta(t)\in H\hat{\otimes} H$ is determined by the formal group law of the Morava $K$-theory. It has finite height $n$ and it is easy to see that $\Phi(\lambda)$, the characteristic polynomial of the Frobenius endomorphism $t\mapsto t^p$ has the form $\lambda^n-p$. 

Recall (cf. for example \cite{Dem}) that there is a contravariant one-to-one correspondence between formal groups of finite height and Dieudonn\'e modules of finite type which are free as $\mathbb{Z}_p$-modules.
We claim that the Dieudonn\'e module corresponding to $H$ has the form $M=R/(V^{n-1}-F)$. Indeed, $M$ clearly has height $n$ and the dimension of the formal group corresponding to $M$ equals the length of the $R$-module $M/VM$ which is equal to one. Now a $1$-dimensional formal group of a given height is determined uniquely by the characteristic polynomial of the Frobenius which is $\lambda^n-p$ in our case. Under the contravariant correspondence between formal groups and Dieudonn\'e modules the Frobenius endomorphism on the formal group side corresponds to the Verschiebung on the Dieudonn\'e module side. Since this characteristic polynomial of the Verschiebung on $M$ is $\lambda^n-p$ we conclude that the Dieudonn\'e module of $H$ is indeed $M$.

Finally, the Dieudonn\'e module ${\mathcal F}(H_\nu)$ is obtained by reducing $M$ modulo $p^\nu$ and we arrive at the following result:
\begin{lemma}
The Dieudonn\'e module ${\mathcal F}(H_\nu)$ is isomorphic to $\mathbb{Z}/p^\nu[F,V]/\sim$ where $\sim$ is  generated by the relations  $V^{n-1}=F$ and $VF=p$.
\end{lemma}\noproof
This lemma allows one to obtain results about the Hopf ring $\Lambda_\boxtimes H_\nu$ via the corresponding Dieudonn\'e
algebra $\Lambda_\boxtimes(M_\nu)$.
\begin{rem}
The ring of coefficients of $K(n)$ is isomorphic to $\F[v_n^{\pm 1}]$ where $v_n$ has degree $2(p^n-1)$. For this version of the Morava $K$-theory its multiplicative structure is determined uniquely and leads to a unique formal group law whose Dieudonn\'e module is as described. Over the field $\F$ there exist many nonisomorphic $1$-dimensional formal groups of height $n$ (all of which become isomorphic after passing to the algebraic closure of $\F$). These formal groups could be realized by the $2$-periodic version of $K(n)$ (which supports many inequivalent product structures). 
\end{rem}

\section{ Structure of the Dieudonn\'e algebra $\Lambda_{\boxtimes}[R/(V^{n-1}-F)]$ }\label{structure}

In this section we investigate the structure of the exterior Dieudonn\'e algebra on the module $M=R/(V^{n-1}-F)$.   It  turns out that with this particular choice of a module this exterior algebra is isomorphic to the conventional exterior algebra on $M$.  Note that the formula for the Dieudonn\'e module of 
$K(n)_*K({\mathbb Z}/p^\nu,q)$ also appears in \cite{SW}

To fix the notation, observe that $M$ is a free ${\mathbb Z}_p$-module with basis $1,V,V^2,\ldots,V^{n-1}=F$.
Denote $V^k$ by $a_{n-k-1}$ where $k=0,\ldots,n-1$. In the basis $a_0,\ldots,a_{n-1}$ the $R$-module structure on $M$ is specified by the formulas:
\begin{align*}Va_0&=pa_{n-1}; Va_i=a_{i-1}, i\geq 1\\
Fa_i&=V^{n-1}a_i, i\geq 0.\end{align*}

Let us consider the (conventional) exterior algebra $\Lambda(M)$ of the free $\mathbb Z_p$-module $M$. Thus, \[\Lambda(M)=\Lambda^0\oplus \Lambda^1\oplus\ldots\oplus\Lambda^n\]
where $\Lambda^q$ is a free $\Z$-module with the basis $a_I=a_{i_1}\wedge\ldots a_{i_q}, 0\leq i_1<
\ldots<i_q<n.$ We will denote the generator of the $1$-dimensional $\Z$-module $\Lambda^0$ by $1$ and that of $\Lambda^n$ - by $a_*=a_0\wedge\ldots\wedge a_{n-1}$.

We will regard $\Lambda^0=\Z$ as the $R$-module $\mathbb I$. Note that $\Lambda^1=M$ has a structure of an $R$-module. Define the structure of an $R$-module on $\Lambda^q$ for $q>1$ inductively by the formulas:
\begin{equation}\label{act1}
Va_I=V(a_{i_1}\wedge\ldots\wedge a_{i_{q-1}})\wedge a_{i_q-1}. \end{equation}
\begin{equation}\label{act2}
 Fa_I=a_{i_1+1}\wedge F(a_{i_2}\wedge \ldots \wedge a_{i_q}).                                                             
\end{equation}
Note that the right hand side of (\ref{act1}) and (\ref{act2}) is well defined since $i_1<n-1$ and $i_q>0$ when $q>1$. Let us check the consistency of the above action. We have:
\begin{align*}
FVa_I&=F(V(a_{i_1}\wedge\ldots\wedge a_{i_{q-1}})\wedge a_{i_q-1})\\
&=(-1)^{q-1}F(a_{i_q-1}\wedge V(a_{i_1}\wedge\ldots\wedge a_{i_{q-1}}))\\
&=(-1)^{q-1}a_{i_q}\wedge FV(a_{i_1}\wedge\ldots\wedge a_{i_{q-1}})\\
&=(-1)^{q-1}pa_{i_q}\wedge(a_{i_1}\wedge\ldots\wedge a_{i_{q-1}})\\
&=pa_{i_1}\wedge\ldots\wedge a_{i_{q}}\\
&=pa_I.
\end{align*}
The condition $VF=p$ is checked similarly. Observe that formulas (\ref{act1}) and (\ref{act2}) imply the following identities:
\begin{equation}\label{act3}
V(a_{i_1}\wedge\ldots\wedge a_{i_q})=Va_{i_1}\wedge\ldots\wedge Va_{i_q}.
\end{equation}
\begin{equation}\label{act4}
F(V(a_{i_1}\wedge\ldots\wedge a_{i_s})\wedge a_{i_{s+1}}\wedge\ldots\wedge a_{i_q})=
a_{i_1}\wedge\ldots\wedge a_{i_s}\wedge F(a_{i_{s+1}}\wedge\ldots\wedge a_{i_q}).
\end{equation}
Note that for these formulas to hold the indices do not necessarily have to be ordered.

The structure of a Dieudonn\'e module on $\Lambda^n$ turns out to be especially simple as the following result demonstrates.
\begin{lemma}
The following formulas hold:
\begin{equation}\label{act5}
 Va_*=(-1)^{n-1}pa_*.
\end{equation}
\begin{equation}\label{act6}
 Fa_*=(-1)^{n-1}a_*.
\end{equation}
\end{lemma}
\begin{proof}
\begin{align*}
F(a_0\wedge\ldots\wedge a_{n-1})=&a_1\wedge F(a_1\wedge\ldots\wedge a_{n-1})\\
=&\ldots=a_1\wedge a_2\ldots\wedge a_{n-1}\wedge Fa_{n-1}\\
=&a_1\wedge a_2\wedge\ldots\wedge a_{n-1}\wedge a_0\\
=&(-1)^{n-1}a_*.
\end{align*}
Next, using formula (\ref{act6}) and the identity $VF=p$ we arrive at formula (\ref{act5}).
\end{proof}
 We will introduce a $\mathbb{Z}_p$-linear scalar product on $\Lambda(M)$ as follows. For $A\in\Lambda^q(L), B\in \Lambda^p(L)$ the product $\langle A,B\rangle=0$ if $p+q\neq n$; the case $p+q=n$ is determined by the formula
\begin{equation}\label{pair}A\wedge B=  \langle A,B\rangle a_*.\end{equation}
\begin{lemma}\label{basic}
The scalar product $\langle,\rangle:\Lambda^q\otimes \Lambda^{n-q}\rightarrow \Z$
satisfies the following property:
\[\langle Va_I,a_J\rangle=(-1)^{n-1}\langle a_I,Fa_J\rangle.\]
\end{lemma}
\begin{proof}We have for $q>0$:
\[Va_I\wedge a_J=\langle Va_I,a_J\rangle a_*.\]
Applying the operator $F$ to the last formula we obtain:
\[F(Va_I\wedge a_J)=\langle Va_I,a_J\rangle Fa_*=(-1)^{n-1}\langle Va_I,a_J\rangle a_*.\]
On the other hand  according to formula (\ref{act4}) we have:
\[F(Va_I\wedge a_J)=a_I\wedge Fa_J=\langle a_I,Fa_J\rangle a_*\]
as required.

Now let $q=0$. For the generator $1\in \Z=\Lambda^0$ we have $V\cdot 1=1$. Therefore
\[\langle V\cdot 1,a_*\rangle=\langle 1,a_*\rangle.\] Using (\ref{act6}) we get
\[\langle 1, Fa_*\rangle=(-1)^{n-1}\langle 1,a_*\rangle.\]
\end{proof}
As a final preparation to our main theorem in this section we will formulate and prove the following general
result.
\begin{lemma}
Let $M,N$ be two $R$-modules which are free as $\Z$-modules. Then the canonical map
\[M\otimes N\rightarrow M\boxtimes N:m\otimes n\mapsto m\circ n=1\otimes m\otimes n\in M\boxtimes N \]
is a monomorphism.
\end{lemma}
\begin{proof}Since $M,N$ have no $p$-torsion it suffices  to prove that the map $M\otimes N\rightarrow M\boxtimes N$ is monomorphic after tensoring with $\mathbb Q$. In other words we have to show that the localization map
\[M\otimes M\rightarrow {\mathbb Q}[V,V^{-1}]\otimes_{{\mathbb Q}[V]}M\otimes M\]
is an inclusion. The latter condition is equivalent to the ${\mathbb Q}[V]$-module $M\otimes N$ having no $V$-torsion. This holds since $M\otimes N$ has no $p$-torsion.   
\end{proof}
\begin{cor}\label{exter}
Let $M$ be a free $\Z$-module. Then the canonical map
\[\Lambda(M)\rightarrow \Lambda_{\boxtimes}(M)\]
is an inclusion.
\end{cor}\noproof
We are now ready to relate $\Lambda(M)$ to the exterior Dieudonn\'e algebra $\Lambda_\boxtimes(M)$.
\begin{theorem}\ 
For any $q=0,1,\ldots, n$ there is an isomorphism of $R$-modules $\phi_q:\Lambda^q(M)\rightarrow\Lambda_{\boxtimes}^q(M)$ given by the formula
\begin{align*}&{\mathbb I}=\Lambda^0\rightarrow\Lambda_\boxtimes^0={\mathbb I}:1\mapsto 1;\\
&\Lambda^q\rightarrow\Lambda_\boxtimes^q:a_{i_1}\wedge\ldots\wedge a_{i_q}\mapsto a_{i_1}\circ\ldots\circ a_{i_q} \textnormal{~~~~~~for~~~~~} q>0.\end{align*}
\end{theorem}
\begin{proof}
A straightforward inspection shows that $\phi$ is indeed a map of $R$-modules. It follows from Corollary \ref{exter} that $\phi$ is an inclusion. It remains to prove that it is an epimorphism. Consider the elements
$a_{i_1}\circ\ldots\circ a_{i_q}\in \Lambda^q_{\boxtimes}$ where $i_1<\ldots<i_q$. We will call such elements \emph{admissible}. Therefore we are reduced to proving that admissible elements additively span $\Lambda^q_{\boxtimes}$. Suppose by induction that this is so and show that admissible elements span  $\Lambda_{\boxtimes}^{q+1}$.

We are assuming that $q>1$ since the cases $q\leq 1$ are obvious.  It follows from the inductive assumption that the collection $\{F^k\otimes a_{i_1}\circ\ldots\circ a_{i_{q+1}}\}$ span $\Lambda_{\boxtimes}^{q+1}$. Here $k=0,1,2\ldots$ and $i_1<\ldots<i_{q+1}$. Since $q>1$ for any given $x=F^k\otimes a_{i_1}\circ\ldots\circ a_{i_{q+1}}$ the elements $a_{i_1},\ldots,a_{i_q}$ are in the image of $V$, namely $a_{i_1}=V(a_{i_1+1}),\ldots, a_{i_q}=V(a_{i_q+1})$.  Using repeatedly the relation (\ref{rel1}) we conclude that $x$ can be rewritten as a multiple of an admissible element. This completes 
the inductive step.
\end{proof}
Now fix a positive integer $\nu$ and consider the Dieudonn\'e module $M_\nu:=M\otimes{\mathbb Z}/p^\nu$ and the corresponding Hopf algebra $H_\nu$. Then the corresponding exterior algebra $\Lambda(M_\nu)\cong \Lambda_{\boxtimes}(M_\nu)$ has a scalar product induced from $\Lambda(M)$ and we arrive at the following result.
\begin{cor}\label{str}\
\begin{enumerate}\item Let $n\geq 1$ be odd, $0\leq q\leq n$. Then the Dieudonn\'e modules $\Lambda^q_{\boxtimes}(M_\nu)$ and $\Lambda^{n-q}_{\boxtimes}(M_\nu)$ as well as Hopf algebras $\Lambda_{\boxtimes}^qH_\nu$ and $\Lambda_{\boxtimes}^{n-q}H_{\nu}$ are dual to each other.\item
Let $n>1$ be even, $0\leq q\leq n$. Then the Dieudonn\'e modules $\Lambda^q_{\boxtimes}(M_\nu)$ and $\Lambda^{n-q}_{\boxtimes}(M_\nu)$ as well as Hopf algebras $\Lambda_{\boxtimes}^qH_\nu$ and $\Lambda_{\boxtimes}^{n-q}H_{\nu}$ are twisted dual to each other. Thus, they become dual in the usual sense after a quadratic extension of $\F$.\end{enumerate}
\end{cor}
\begin{proof}
It follows from Lemma \ref{basic} that the scalar product on $\Lambda_{\boxtimes}(M)$  induces a bilinear pairing between $\Lambda^{q}_{\boxtimes}(M_\nu)$ and $\Lambda^{n-q}_{\boxtimes}(M_\nu)$ for $n$ odd. Similarly for $n$ even this pairing is twisted bilinear. It is clearly nondegenerate and the result follows. 
\end{proof}
We conclude this section with the useful observation that the $R$-module $\Lambda^q(M)\cong \Lambda^q_{\boxtimes}(M)$ is actually a module over $R/(V^{n-q}-F^q)$.
\begin{prop}\label{module}
In the $R$-module $\Lambda^q(M)$ the following relation holds: $V^{n-q}=F^q$.
\end{prop}
\begin{proof}
Note that since $VF=p$ and $\Lambda^q(M)$ is a free $\Z$-module the operator $V$ is monomorphic on $\Lambda^q(M)$. Therefore it suffices to prove the relation $V^q\cdot V^{n-q}=V^q\cdot F^q$ or $V^n=p^q$.
For a generator $a_{i_1}\wedge\ldots\wedge a_{i_q}$ we have according to formula (\ref{act3})
\begin{align*}V^n(a_{i_1}\wedge\ldots\wedge a_{i_q})=&V^na_{i_1}\wedge\ldots\wedge V^na_{i_q}\\
=&pa_{i_1}\wedge\ldots\wedge pa_{i_q}\\
=&p^qa_{i_1}\wedge\ldots\wedge a_{i_q}\end{align*}
as required.
\end{proof}
\begin{rem}
We want to stress that the isomorphism of $\Lambda_\boxtimes(M)$ with the usual exterior algebra  $\Lambda(M)$ relies essentially on the peculiar properties of the module $M=R/(V^{n-1}-F)$. Calculations show that bilinear products as well as exterior Dieudonn\'e algebras of cyclic $R$-modules which are free of finite rank over $\Z$ usually contain $p$-torsion and are often no longer finitely generated over $\Z$. 
\end{rem}
\section{Decomposition of $\Lambda_{\boxtimes}[R/(V^{n-1}-F)]$ up to isogeny over $\bar{\mathbb F}_p$}

Here $\mathbf k$ will be an algebraic extension of $\F$, in fact we will be most interested in the case when 
${\mathbf k}=\bar{\mathbb F}_p$, the algebraic closure of $\F$. Recall the definition of the Dieudonn\'e ring $R(\mathbf k)$ and $R(\mathbf k)$-module $N(\mathbf k)$ for an $R$-module $N$ from Section \ref{twisted}. We will denote by $\widetilde{W}(\mathbf k)$ the field of fractions of $W(\mathbf k)$. The ring $\widetilde{W}(\mathbf k)\otimes_{W(\mathbf k)} R\cong {\mathbb Q}\otimes R(\mathbf k)$ will be denoted by $\tilde{R}(\mathbf k)$ and $\tilde{N}(\mathbf k)$ will stand for the $\tilde{R}(\mathbf k)$-module $\widetilde{W}(\mathbf k)\otimes_{W(\mathbf k)} N\cong {\mathbb Q}\otimes N(\mathbf k)$.
\begin{defi}
An $F$-space is a $\tilde{R}(\mathbf k)$-module which is finite-dimensional as a 
$\widetilde{W}(\mathbf k)$-vector space. 
$F$-spaces form a category whose morphisms are simply the morphisms of $\tilde{R}(\mathbf k)$-modules.
\end{defi}                                                                       
The $R(\mathbf k)$-modules we will be interested in here will be torsion-free and finite rank over ${W}(\mathbf k)$. Any such module $N$ determines a lattice in the $F$-space $\tilde{N}=\mathbb{Q}\otimes N$. 
\begin{defi}
Let $N_1$ and $N_2$ be two torsion-free $R(\mathbf k)$ modules of finite rank over ${W}(\mathbf k)$. We say that $N_1$ and $N_2$ are \emph{isogenous} if the corresponding $F$-spaces $\tilde{N}_1$ and $\tilde{N}_2$ are isomorphic.
\end{defi}
Now let us introduce the Dieudonn\'e modules $R_{n,q}=R/(V^{n-q}-F^q)$.
Observe that  \[\tilde{R}_{n,q}=\widetilde{W}({\mathbf k})[T]/(p^{n-q}-T^n)\]
and $F$ acts as a multiplication by $T$.
It is interesting to note that for $p, q$ coprime and $\mathbf k$ an algebraic extension of $\F$ of degree $n$ the module $\tilde{R}_{n,q}(\mathbf k)$ is a division algebra over $\mathbb Q_p$ with Hasse invariant $q/n$, see. e.g. \cite{Pierce}, Chapter 17. The module  $R_{n,q}(\mathbf k)$ is in this case a maximal order in $\tilde{R}_{n,q}(\mathbf k)$.

From now on we assume that $\mathbf k =\F$. 
It is known that the category of $F$-spaces
is semisimple and all simple objects are of the form $\tilde{R}_{n,q}$ where $n,q$ are nonnegative relatively prime integers and $n>0$. Recall that we denoted by $M$ the $R$-module $R_{n,1}=R/(V^{n-1}-F)$. Let $(n,q)$ be the greatest common factor of $n$ and $q$ and set $n_0:=\frac{n}{(n,q)}$, $q_0:=\frac{q}{(n,q)}$ We have the following result.
\begin{theorem}\label{isogeny}
The ${R}(\mathbf k)$-module $\Lambda^q(M)$ is isogenous to the direct sum of $\frac{1}{n_0}{n\choose q}$ copies of the $R(\mathbf k)$-module $R_{n_0,q_0}$.
\end{theorem}
\begin{proof}
Observe that the operator $V$ permutes (up to multiplication by a scalar factor) the basis vectors $a_I$ in the  $\tilde{R}$-module ${\mathbb Q}\otimes\Lambda^q(M)$.  Since $V^n=p^q$ in ${\mathbb Q}\otimes\Lambda^q(M)$ by Proposition \ref{module} we see that this module decomposes into a direct sum of submodules isomorphic to quotients of $\tilde{R}_{n,q}$.  Similarly ${\mathbb Q}\otimes \Lambda^q(M)(\mathbf k)$ decomposes into a direct sum of quotients of  $\tilde{R}_{n,q}(\mathbf k)$. Therefore it suffices to show that $\tilde{R}_{n,q}(\mathbf k)$ is isomorphic to the direct sum of copies of $\tilde{R}_{n_0,q_0}(\mathbf k)$. We claim that if there exists a nontrivial map of $R$-modules \begin{equation}\label{ma}\tilde{R}_{r,s}(\mathbf k)\rightarrow\tilde{R}_{n,q}(\mathbf k)\end{equation} then $\frac{r}{s}=\frac{n}{q}$.  This claim will clearly imply what we need.

The proof of the claim is similar to Proposition D in \cite{Dem}, page 79.  A map (\ref{ma}) is equivalent to having an element $x\in \tilde{R}_{n,q}(\mathbf k)$ for which $F^rx=p^{r-s}x$. Note that $\tilde{R}_{n,q}(\mathbf k)$ has a basis $f_j$ such that if $x=\sum b_jf_j$ then \[F^{n}x=\sum\sigma^n(b_j)p^{n-q}f_j\] and so 
\[F^{rn}x=\sum\sigma^{rn}(b_j)p^{r(n-q)}f_j.\]
On the other hand if $F^rx=p^{r-s}x$ then 
\begin{align*}F^{rn}x=&p^{(r-s)n}x\\
=&\sum b_jp^{(r-s)n}f_j.\end{align*}
It follows that $\sigma^n(b_j)p^{(r-s)n}=b_jp^{r(n-q)}$ and since the Frobenius $\sigma$ preserves $p$-adic valuation on $\widetilde{W}(\mathbf k)$ we conclude that $(r-s)n=r(n-q)$. It follows that $\frac{r}{s}=\frac{n}{q}$
as required.
\end{proof}

\section{ $p$-divisible groups associated with $K(n)^*K(\mathbb Z, q)$}
\subsection{Basic definitions} We start by recalling some standard definitions and facts from the theory of $p$-divisible groups referring the reader to \cite{Dem}, \cite{tate} or \cite{shatz} for details. Here $\mathbf k$ is an algebraic extension of $\F$.
\begin{defi} A $p$-divisible group of height $h$ over a field $\F$
is a sequence $G=(H_\nu, i_\nu)$, $\nu=0,1,2,\ldots$ of Hopf algebras over ${\mathbf k}$ with $\dim{H_\nu}=\nu h$ and $i_\nu:H_{\nu+1}\rightarrow H_\nu$ is a Hopf algebra homomorphism such that for each $\nu$ the sequence
\[\xymatrix{H_{\nu+1}\ar^{[p^\nu]}[r]&H_{\nu+1}\ar^{i_\nu}[r]&H_\nu\ar[r]&0}\]
is exact in $\A$.\end{defi}

Since $[p^{\nu+1}]=[p]\cdot[p^\nu]$ is the trivial endomorphism of $H_{\nu+1}$ we conclude that there exists a map $j_{\nu}:H_\nu\lra H_{\nu+1}$ making commutative the following diagram
\[\xymatrix{H_{\nu+1}\\H_{\nu+1}\ar^{[p]}[u]\ar^{i_\nu}[r]&H_{\nu}\ar_{j_\nu}[ul]}\]

The topological Hopf algebra $H=\inlim{}H_\nu$ represents a formal group from which the sequence  $(H_\nu, i_\nu)$ could be recovered by setting $H_\nu:=\coker\{[p^\nu]:H\rightarrow H\}$. We will use the term `$p$-divisible group' also for the corresponding formal group.   The dimension of the $p$-divisible group $G$ is the Krull dimension of $H$. 

Next, $H$ will be isomorphic to a ring of formal power series (one say that in this case the corresponding formal group is \emph{smooth}) if and only if each $H_\nu$ is a local ring.

There is a suitable version of duality for $p$-divisible groups. As far as we know twisted duality has not been considered before.

\begin{defi}If $G=(H_\nu, i_\nu)$ is a $p$-divisible group then its \emph{Serre dual} $p$-divisible group
is defined as $DG=(H_\nu^*, j_\nu^*)$. It has the same height as $(H_\nu, i_\nu)$. Its \emph{twisted Serre dual} is the collection $D^\prime G=(D^\prime H_\nu, D^\prime j_\nu)$ where $D^\prime$ is the functor of twisted duality on the category of $p$-torsion Hopf algebras.\end{defi}
\begin{rem}
Of course, the notions of the Serre dual and twisted Serre dual $p$-divisible group coincide when ${\mathbf k}$ contains a quadratic extension of $\F$.
\begin{example}\
\begin{enumerate}\item A one-dimensional $p$-divisible group of height $1$ is represented by  the topological Hopf algebra $H={\mathbf k}[[x]]$ with $\Delta(x)=1\otimes x+x\otimes 1+x\otimes x$, the so-called multiplicative formal group. This $p$-divisible group is smooth. Its Serre dual $p$-divisible group is represented by the Hopf algebra ${\mathbf k}[{\mathbb Z}]^*$, the dual group ring of the infinite cyclic group. The  dimension of this $p$-divisible group is zero and it is not smooth.
\item Let $X$ be an abelian scheme over ${\mathbf k}$. Then the kernel of the multiplication by $p^\nu$ on $X$ is a finite group scheme which is represented by a Hopf algebra whose dimension is a power of $p$. The resulting inverse system of Hopf algebras constitutes a $p$-divisible group. \end{enumerate}
\end{example}  
\end{rem}
To conclude our review of the background material we note that the height of a $p$-divisible group is equal to the sum of its dimension and the dimension of its dual:
\[\textnormal{height}(G)=\textnormal{height}(DG)=\dim(G)+\dim(DG).\]
\subsection{Dieudonn\'e correspondence} We now briefly review the Dieudonn\'e correspondence adapted to $p$-divisible groups following \cite{Dem}. Let $G=(H_\nu, i_\nu)$ be a $p$-divisible group. Its Dieudonn\'e module
is defined as ${\mathcal F}(G)=\inlim{}{\mathcal F}(H_\nu)$. Then we have a theorem:
\begin{theorem}
The correspondence $G\mapsto {\mathcal F}(G)$ induces an equivalence between the categories of $p$-divisible groups over $\mathbf k$ and $R(\mathbf k)$-modules which are finite rank free $W({\mathbf k})$-modules. The height of $G$ is equal to the dimension of ${\mathcal F}(G)$.  If the Dieudonn\'e module corresponding to $G$ is actually a module over $\hat{R}(\mathbf k):=\inlim{n}R(\mathbf k)/V^n$ then $G$ is smooth. The dimension of $G$ equals $\dim_{\mathbf k}[({\mathcal F}(G)/V{\mathcal F}(G)]$.
\end{theorem}
\begin{example}
Let $R_{n,q}=R({\mathbf k})/(V^{n-q}-F^q)$. The corresponding $p$-divisible group $G$ has height $n$. It is smooth of dimension $q$ for $q=1,2\ldots, n$.  The Serre dual to $G$ has Dieudonn\'e module $R_{n,n-q}=R({\mathbf k})/(V^{q}-F^{n-q})$.
\end{example}

For any $p$-divisible group $G$ one can define the $F$-space 
$ {\mathcal F}(G)\otimes \mathbb Q$. 
Then two $p$-divisible groups are isogenous (i.e. there exists a monomorphism between their representing Hopf algebras with a finite-dimensional cokernel) if and only if the corresponding $F$-spaces are isomorphic.

\subsection{Main theorems}
Now let $H_\nu(q):=K(n)^*K({\mathbb Z}/p^\nu, q)$. Also denote $H_\nu^*(q)=K(n)_*K({\mathbb Z}/p^\nu, q)$ by
$H^\nu(q)$. The inclusion \[\mathbb{Z}/p^\nu\lra  \mathbb{Z}/p^{\nu+1}\]
induces a map of spaces \[K({\mathbb Z}/p^\nu, q) \lra K({\mathbb Z}/p^{\nu+1}, q)\] which in turn gives rise
to a map of Hopf algebras  $i_\nu:H_{\nu+1}(q)\rightarrow H_\nu(q)$.

We can formulate now our main result.
\begin{theorem}\label{main}\
\begin{enumerate}\item The sequence  $(H_\nu(q), i_\nu)$ forms a $p$-divisible group of height $n\choose q$.
 \item For $q=1,2,\ldots,n-1$ the $p$-divisible group $(H_\nu(q), i_\nu)$ is smooth. The corresponding formal group is represented by a formal power series ring on $n-1\choose q-1$ variables and could be identified with
$K(n)^*K(\mathbb{Z}, q+1)$.\item If $n$ is odd then there is an isomorphism \[H^\nu(n)\cong D(\nu):=\F[\mathbb{Z}/p^\nu]^*.\] If $n$ is even then $H^\nu(n)\cong D^\prime(\nu)$ where $D^\prime(\nu)$ is the twisted dualizing Hopf algebra.
\item \begin{enumerate}\item Let $n$ be odd and $0<q<n$.
 Then the $\circ$-pairing \[H^\nu(q)]\boxtimes H^{\nu}(n-q)\lra H^\nu(n)\] induces an isomorphism of the formal group of $K(n)^*K(\mathbb{Z}, q+1)$ with the Serre dual of the formal group of
$K(n)^*K(\mathbb{Z},n-q+1 )$. In particular the Hopf algebras 
$H^\nu(q)$ and $H_\nu(n-q)$ are isomorphic. \item  For $n$ even the above statement holds provided one replaces `dual' with `twisted dual'.   \end{enumerate}
 \end{enumerate}
\end{theorem}
\begin{proof}
We saw in section \ref{structure} that the Dieudonn\'e modules corresponding to $K({\mathbb Z}/p^\nu, q)$ 
are obtained by reducing modulo $p^\nu$ from the modules $\Lambda_{\boxtimes}^q(M)$ where $M=R/(V^{n-1}-F)$. 
These modules torsion-free and their rank over $\Z$ is equal to $\dim\Lambda^q(M)={n\choose q}$. This proves part (1). 

Further note that since the identity $V^{n-q}=F^q$ holds in $\Lambda_{\boxtimes}(M)\cong \Lambda(M)$ the corresponding $p$-divisible group is smooth for $q<n$. From the formula (\ref{act3}) for the action of $V$ on $\Lambda_{\boxtimes}(M)$ we deduce that the image of $V$ is spanned modulo $p$ by the monomials $a_{i_1}\wedge \ldots\wedge a_{i_q}$ where $i_q\neq n-1$. It follows that the monomials $a_{i_1}\wedge \ldots\wedge a_{i_{q-1}}\wedge a_{n-1}\mod V$ form a basis in $\Lambda^q(M)/V\Lambda^q(M)$ and so its dimension over $\F$ equals $n-1\choose q-1$. Therefore the dimension of the formal group in question is as claimed.

Next, observe that there is a homotopy equivalence of spaces
\[K({\mathbb Q/\mathbb Z},q)\simeq \hocolim K({\mathbb Z/p^\nu},q)\]
where the homotopy colimit is taken over all $\nu$ and prime numbers $p$.
This homotopy equivalence induces an isomorphism of Hopf algebras
\[\inlim{\nu}K(n)^*K({\mathbb Z}/p^\nu,q)\cong K(n)^*K({\mathbb Q/\mathbb Z},q).\]
The obvious map $K({\mathbb Q/\mathbb Z},q))\rightarrow K({\mathbb Z},q+1)$ induces an isomorphism on $K(n)$-theory. Therefore the $p$-divisible group $H_\nu=K(n)^*K({\mathbb Z}/p^\nu,q)$ is indeed representable by the Hopf algebra  $K(n)^*K({\mathbb Z},q+1)$ as claimed. This proves part (2). The remaining claims (3) and (4) follow from Corollary \ref{str}.\end{proof}
\begin{rem}
It is curious to note that the Serre duality between $K(n)^*K(\mathbb{Z}, q+1)$ and $K(n)^*K(\mathbb{Z},n-q+1 )$ breaks down for $q=0$. Indeed, $K(\mathbb{Z},1)$ is simply the circle $S^1$ and so $K(n)^*K(\mathbb{Z}, 1)$ cannot give rise to a formal group. However the duality between $K(n)_*K(\mathbb{Z}/p^\nu, 0)$ and $K(n)_*K(\mathbb{Z}/p^\nu, n)$ continues to hold. The point is that $K(n)^*K(\mathbb{Z}, q+1)$ can no longer be related to  $K(n)_*K(\mathbb{Z}/p^\nu, q)$ for $q=0$.
\end{rem}
\begin{rem}
Our results are formulated under the assumption that $\mathbf k$ has characteristic $p\neq 2$.  Note, however, that an appropriate version of Theorem \ref{main} holds for $p=2$ as well. The point is that although in this case $K(n)$ is not a commutative ring spectrum, $K(n)_*K({\mathbb Z}/2^\nu,q)$ is still a bicommutative Hopf algebra. See \cite{JW} for details. In the exterior Dieudonn\'e algebra we should impose the additional condition that the circle product of two equal elements is zero. The notion of twisted duality is extraneous; the twisted dualizing object is the same as the untwisted one.
\end{rem}
The identification of the isogeny class of $K(n)^*K(\mathbb Z, q+1)$ is an immediate consequence of Theorem \ref{isogeny}. Recall that for two nonnegative integers $q,n$ the pait $n_0, q_0$ is specified by the condition that $\frac{q}{n}=\frac{q_0}{n_0}$ and that $q_0,n_0$ be coprime.  Thus, we obtain the following result.
\begin{theorem}  
For $0<q<n$ the formal group of $K(n)^*K(\mathbb Z, q+1)$ is isogenous to the product of $ \frac{1}{n_0}{n\choose q}$ copies of the $p$-divisible group corresponding to the Dieudonn\'e module $R_{n_0,q_0}$.
\end{theorem}\noproof
\begin{rem}
One might wonder whether the formal groups corresponding to $K(n)^*K(\mathbb Z, q+1)$ are \emph{algebraicizable}, i.e. whether there exists abelian schemes of which they are formal completions. Since
abelian schemes are always isogenous with their dual the same is true for their completions. This is the so-called Manin symmetry condition, see \cite{man} or \cite{oort}. It follows that if $n$ is odd or if $n$ is even but $q\neq \frac{n}{2}$ the formal group of $K(n)^*K(\mathbb Z, q+1)$ cannot be algebraicized. Note that for $n$ even the formal group of $K(n)^*K(\mathbb Z, \frac{n}{2}+1)$ is \emph{supersingular}, i.e. it is isogenous (over $\bar{\mathbb F}_p$) to the product of copies of one-dimensional formal group of height 2.
\end{rem}

\end{document}